\documentclass[12pt]{amsart}

\def\CA{{\mathcal A}}
\def\CH{{\mathcal H}}

\def\acts{\triangleright}

\def\id{\hbox{id}}

\def\ul#1{\underline{#1}}
\def\lm#1{\lambda^{#1}}
\theoremstyle{remark}
\newtheorem{remark}[]{Remark}
\newtheorem{lemma}[]{Lemma}
\newtheorem{definition}[]{Definition}
\def\ts{\otimes}
\def\cop{\triangle}

\def\acts{\triangleright}

\def\endproof{\vrule height 0.5em depth 0.2em width 0.5em}
\title{Twists and spectral triples for isospectral deformations.}
\author{Andrzej Sitarz}
\email{Andrzej.Sitarz@th.u-psud.fr}
\thanks{${}^\dagger$\ Supported by Marie Curie Fellowship.}
\address{Laboratoire de Physique Theorique, Universit\'e Paris-Sud, 
Bat. 210, 91405  ORSAY Cedex, France}
\address{Institute of Physics, Jagiellonian University, Reymonta 4, 
30-059 Krak\'ow, Poland}
\begin{document}
\begin{abstract}
We construct explicitly the symmetries of the isospectral 
deformations as twists of Lie algebras and demonstrate 
that they are isometries of the deformed spectral triples.
\end{abstract}
\maketitle
\section{Introduction}
The isospectral deformations, which have been introduced 
by Connes and Landi \cite{CoLa}, with the examples of the 
noncommutative 3 and 4 spheres, have attracted recently 
much attention. 

The deformation itself, which, looking only at the algebra is 
of the "star deformation" (Moyal) type, has been discussed
in other contexts and the construction of similar, physically
motivated examples can be found in recent works by Kulish 
and Mudrov \cite{Kul} as well as in Paschke Ph.D. thesis \cite{Pas}. 

On the mathematical side (though also with a physical motivation) 
it could be traced back to sixties \cite{GLS} or, more recently, 
to the works of Rieffel \cite{Rie1,Rie2,Rie3}, Dubois-Violette 
\cite{MDV} and Mourre \cite{Mor}.

However, only after \cite{CoLa} it became apparent that these
kind of deformation allows to describe a real noncommutative
spin manifold in the sense of real spectral triples.  Moreover,
their geometry is distinct from the most of other known deformations 
allowing for definition of $C^\infty$ elements, no dimension drop 
and the transfer of most geometric features from the 
commutative case \cite{CMDV}.

In this paper we discuss the notion of symmetries of the
spectral triples for isospectral deformations. The existence 
of quantum group symmetries for isospectral deformations 
has been first noted by Connes and Dubois-Violette \cite{CMDV}. 
Here, for the purpose of making connections with the work 
on symmetries of spectral triples \cite{PaSi} we have 
independently developed a dual approach (using Lie 
algebras and not Lie groups). We demonstrate that the 
twist by a Cartan subalgebra provides a deformed 
symmetry of the {\em isospectral deformation} and that it 
is an isometry of the deformed manifold.

\section{Symmetries of the isospectral deformations}
The isospectral deformation has appeared in the construction
of the examples of noncommutative 4 spheres, which have 
the same instanton bundle as the "classical" sphere. One of 
the  crucial questions posed by the construction was whether 
the constructed {\em spheres} still are symmetric, i.e.~whether 
the natural $SO(5)$ symmetries (or, respectively, $SO(4)$ for 
the 3-sphere) are preserved (in the form of deformed Hopf 
algebras) and whether the constructed spectral triples are 
{\em symmetric} in the sense of \cite{PaSi}. 

Here, we shall present the symmetry in terms of the deformation 
of the universal enveloping algebra acting on the deformed algebra. 
We shall use both a general approach as well as the description in 
terms of generators and relations, the latter to make connections 
with the known specific examples. This enabled us to make
connection \cite{Sit} with the work of Rieffel \cite{Rie1,Rie2} and 
to present the isospectral deformation (on the algebraic level) as 
a special case of the strict deformation quantization as 
described elsewhere \cite{Rie3}.

Let $H$ be a Hopf algebra (symmetry algebra), which contains two 
independent, mutually commuting generators of $U(1)$ symmetries,  
i.e. operators $h_1, h_2$ such that: 
$$ [ h_1, h_2] = 0, \;\;\;\; \cop h_i = 1 \ts h_i + h_i \ts 1. $$
Consider an algebra $\CA$ on which $H$ acts 
from the left. The generators $h_i$ act on $\CA$ as generators 
of mutually independent $U(1)$ symmetries. 

\begin{remark}
With respect to the action of $h_1,h_2$ we can select out 
elements in $\CA$ of degree $(n_1,n_2)$, we say that
$t \in \CA$ is of degree $(n_1,n_2)$ iff:
\begin{equation}
\begin{array}{ll}
h_1 \acts t = n_1 t, & h_2 \acts t = n_2 t.
\end{array}
\end{equation}
\end{remark}

\begin{remark}
The product in the algebra $\CA$ can be deformed, first on
elements of given degree, and then extending the deformation
by linearity:
\begin{equation}
a * b = ab \lambda^{n_1^a n_2^b}, \label{deform}
\end{equation}
where $\lambda$ is a complex number such that $|\lambda|=1$. 

This is the crucial step of the construction of the isospectral 
deformation as described in \cite{CoLa}, for future reference 
we shall denote the deformed algebra by $\CA_\lambda$.
\end{remark}

This deformation, can be extended, in the case 
of a differential manifold and a Lie algebra acting 
on the differential functions,  to all $C^\infty$ functions 
\cite{CoLa}.

It would be useful to introduce a {\em quantization map}:
$$ \CA \ni a  \mapsto \ul{a} \in \CA_\lambda, $$
so that the Eq.\ref{deform} could be rewritten as:
\begin{equation}
\ul{a} * \ul{b} = \ul{ab} \lm{n_1^a n_2^b}, 
\label{deform2}
\end{equation}

Before we proceed with the deformation of the 
symmetry and the generalization of the above
deformation, let us observe that the star structure 
of the algebra $\CA_\lambda$ could be also deformed:
\begin{equation}
(\ul{a})^\ast = \lm{n_1^a n_2^a} \ul{(a^\ast)}, \label{star}
\end{equation}

\subsection{Example: Lie algebras of rank $n \geq 2$}

Let $U(\mathfrak g)$ be the universal enveloping algebra of 
a Lie algebra $\mathfrak g$ with the Cartan matrix $a_{ij}$, 
and the set of generators in the Chevalley basis: $h_i$, 
which commute with each other (we assume that the rank 
of the Lie algebra is at least 2). 
\begin{equation}
[h_i, h_j] = 0, \label{ug1}
\end{equation}
and  $x_k^\pm$, which satisfy:
\begin{eqnarray}
&& [h_i, x_j^\pm] = \pm a_{ij} x_j^\pm. \label{ug2} \\
&& [x_i^+, x_j^-] = \delta_{ij} h_i, \label{ug3}
\end{eqnarray}
together with the Serre relations:
\begin{equation}
\sum_{k=0}^{1-a_{ij}} \binom{1-a_{ij}}{k} 
(x_i^+)^k x^+_j (x_i^+)^{1 - a_{ij}-k} = 0.
\label{serre}
\end{equation}

Then, choosing, for instance, $h_1,h_2$ as the generators 
discussed earlier we have the desired situation.

\section{The twisted symmetry $H_\lambda$}
We shall construct here the deformation of the symmetry 
algebra $H_\lambda$, which is the twist (cocycle deformation) 
of $H$. Furthermore, we shall 
verify in the next section that this will be the 
symmetry of the deformed algebra $\CA_\lambda$.
It will have the same algebra structure and the same action 
of the generators on the elements of $\CA_\lambda$, with the 
coalgebra and antipode modified by  a twist.  For a general 
theory, examples and details of twists see \cite{CP} and \cite{Maj}.
We shall use capital letters to denote the elements of the 
deformed algebra, then, if no misinterpretation is possible 
we shall write $\cop$ instead of $\cop_\lambda$.
 
\begin{definition}
Let us define $H_\Psi$ as an algebra isomorphic to $H$, 
however, with a twisted coproduct:
\begin{equation}
\cop_\lambda T = \Psi \cop t \Psi^{-1}, \label{cotwi-0}
 \end{equation}
 where $\Psi$ is an invertible element of $H \ts H$, which
 satisfies a cocycle condition:
\begin{eqnarray}
&& \Psi_{12} ( \cop \ts \id) \Psi = \Psi_{23} (\id \ts \cop) \Psi, \label{twi-1} \\
&& (\epsilon \ts \id) \Psi =1 = (\id \ts \epsilon) \Psi.
\end{eqnarray}
\end{definition} 
In the particular example associated with the isospectral
deformations we shall take an element $\Psi_c$ associated 
with the Cartan subalgebra generated by $h_1,h_2$:
$$H \ts H \ni \Psi_c = \lm{-H_1 \ts H_2},$$
and we shall call this particular twist $H_\lambda $. 
The deformed coproduct in $H_\lambda$ becomes:
\begin{equation}
\cop_\lambda T = \lm{-H_1 \ts H_2} (\cop t) \lm{H_1 \ts H_2}, \label{cotwi}
\end{equation}

To see that $\Psi_c$ satisfies the cocycle condition 
we shall verify it explicitly:

$$
\begin{array}{l}
\Psi_{12} ( \cop \ts \id) \Psi = \lm{-H_1 \ts H_2 \ts 1} 
\left( \cop \ts \id \right) \lm{-H_1 \ts H_2} = \\ 
\phantom{xxx} = \lm{-H_1 \ts H_2 \ts 1} \lm{ - (1 \ts H_1 + H_1 \ts 1) \ts H_2 } = \\
\phantom{xxx} = \lm{ - H_1 \ts H_2 \ts 1 - 1 \ts H_1 \ts H_2 - H_1 \ts 1 \ts H_2}.
\end{array}
$$
on the other hand:
$$
\begin{array}{l}
\Psi_{23} ( \id \ts \cop) \Psi = \lm{-1 \ts H_1 \ts H_2} 
\left( \id \ts \cop \right) \lm{-H_1 \ts H_2} = \\ 
\phantom{xxx} = \lm{- 1 \ts H_1 \ts H_2} \lm{ -H_1 \ts  (1 \ts H_2 + H_2 \ts 1)} = \\
\phantom{xxx} = \lm{ - H_1 \ts H_2 \ts 1 - 1 \ts H_1 \ts H_2 - H_1 \ts 1 \ts H_2}.
\end{array}
$$
which proves (\ref{twi-1}). \endproof

The counit does not change:
\begin{equation}
\epsilon(T) = \epsilon(t),
\end{equation}
whereas the antipode is twisted by an element $U$:
\begin{equation}
U = \Psi_{(1)} S( \Psi_{(2)} ), \label{utwi}
\end{equation}
\begin{equation}
S(T) = U S(t) U^{-1}. \label{twianti}
\end{equation}

\subsection{Twist of $U_\lambda({\mathfrak g})$}

In the particular example of the Lie algebra 
${\mathfrak g}$, we have an algebra with the same relations 
as (\ref{ug1}-\ref{serre}) . For simplicity we shall 
call $a_{1i}=\alpha_i$ and $a_{2i}= \beta_i$.  

Then, using (\ref{cotwi}) we might calculate explicitly:
\begin{equation}
\begin{array}{l}
 \cop {H_i} = {H_i} \ts 1 + 1  \ts {H_i}, \\
 \cop {X_i^\pm} = \lm{-H_1 \ts H_2} \left( {X_i^\pm} \ts 1
+ 1 \ts {X_i^\pm} \right)  \lm{-H_1 \ts H_2} = \\
\phantom{xxx} = {X_i^\pm} \ts \lm{\mp \alpha_i H_2}
+ \lm{\mp \beta_i H_1} \ts {X_i^\pm}, 
\end{array} \label{co2}
\end{equation}

The obtained object is a triangular Hopf algebra 
with the universal ${\mathcal R}$-matrix:
$$ {\mathcal R} = \lm{H_2 \ts H_1 - H_1 \ts H_2}. $$

For completeness we give here the counit and 
the antipode calculated using (\ref{twianti}):

\begin{equation}
\begin{array}{lp{2cm}l}
\epsilon(H_i) =0, & & \epsilon(X_i^\pm) =0, \\
S H_i = - H_i & & S X_i^\pm = - \lm{\pm \beta_i H_1} 
X_i^\pm \lm{\pm \alpha_i H_2},
\end{array}
\end{equation}
\section{Twisted symmetry of the algebra $\CA_\lambda$}

We shall demonstrate now that the twisted Hopf algebra defined 
in the previous section is the symmetry algebra of the deformation 
$\CA_\lambda$. 

Even more, to see it, we shall use the generalization 
of the deformation of $\CA$ by an arbitrary twist $\Psi$.

\begin{definition}
Let us define a deformed product for elements of $\CA$ 
(we shall denote the deformed algebra $\CA_\Psi$):
If $\ul{a},\ul{b} \in \CA_\Psi$ then:
\begin{equation}
\ul{m} \left( \ul{a} \ts \ul{b} \right) = 
 \ul{m \left( \Psi^{-1} \acts (a \ts b) \right)}, \label{prod}
\end{equation}

where $m$ is the multiplication map 
$$m:\CA \ts \CA \ni a \ts b \to ab \in \CA. $$
and $\ul{m}$ is the multiplication in $\CA_\Psi$:
$$\ul{m}:\CA_\Psi \ts \CA_\Psi \ni \ul{a} \ts \ul{b} \to 
\ul{a}* \ul{b} \in \CA_\Psi. $$
\end{definition}
We shall verify explicitly this defines an associative product:
$$
\begin{array}{l}
\ul{a} * (\ul{b} * \ul{c})  = \ul{ m  \left( \Psi^{-1} \acts (a \ts (\ul{b} * \ul{c} \right) }= \\
\phantom{xxx} = \ul{ ( (\id \ts \cop) \Psi^{-1}) \Psi^{-1}_{23} \acts (a \ts b \ts c) } = \ldots 
\end{array} 
$$
and using the cocycle condition for $\Psi$;
$$
\begin{array}{l}
\phantom{xx} \ldots = \ldots  \ul{ ( (\cop \ts \id) \Psi^{-1} ) \Psi^{-1}_{12} 
 \acts (a \ts b \ts c) }= \\
\phantom{xxx} =( \ul{a} * \ul{b})  * \ul{c}, 
\end{array} $$
\endproof

\begin{definition}
The action of the twisted symmetry  $H_\Psi$
on $\CA_\Psi$ could be defined as the twisting 
of the action of $H$ on $\CA$. If $\ul{a} \in \CA_\Psi$,
and $T \in H_\Psi$ then:

\begin{equation}
T \acts \ul{a} = \ul{(t \acts a)}, \label{action}
\end{equation}
where $T$ is deformed $t$.
\end{definition}

To verify that this is an action of the Hopf algebra on 
$\CA_\Psi$ we have to verify first the compatibility with the product.  

We shall present here a general proof, which we shall repeat later 
in the specific case of $U_\lambda({\mathfrak g})$, using generators 
and homogeneous elements. The general proof uses the twisting relation 
(\ref{cotwi}).

\begin{proof}
We have, for arbitrary $T \in H_\Psi$ (which we
identify with $t \in H$):
\begin{equation}
\begin{array}{l}
T \acts \left( \ul{a} * \ul{b} \right) = \\
\phantom{xxx} = \ul{ t \acts m \left( 
\Psi^{-1}  \acts (a \otimes b) \right) }= \\
\phantom{xxx} \ul{ m \left( \left( ( t_{(1)} \ts t_{(2)} ) 
\Psi^{-1} \right) \acts (a \otimes b) \right)},
\end{array}
\end{equation}
where we have used the form of the product (\ref{prod}) the
definition of the action (\ref{action}) and the undeformed 
coproduct of $t$. On the other hand we have:
\begin{equation}
\begin{array}{l}
T \acts \left( \ul{a} * \ul{b} \right) = \\
\phantom{xxx} = ( T_{(1)} \acts \ul{a} ) * ( T_{(2)} \acts \ul{b} )  =\\
\phantom{xxx} = \ul{m} \left( \left(
\Psi ( t_{(1)} \ts t_{(2)} ) \Psi^{-1} \right) \acts
( \ul{a} \ts \ul{b} ) \right)= \\
\phantom{xxx} =
m \ul{ \left( \left( \Psi^{-1} \Psi 
( t_{(1)} \ts t_{(2)} ) \Psi^{-1} \right) \acts ( a \ts b) \right] } = \\
\phantom{xxx} =
m \ul{ \left( ( t_{(1)} \ts t_{(2)} ) \Psi^{-1} \right) \acts ( a \ts b) },
\end{array}
\end{equation}
which ends the proof.
\end{proof}

\begin{remark}
In particular, $H_\lambda$ is the symmetry of the deformed 
$\CA_\lambda$ algebra.
\end{remark}

\section{The star structure}

Suppose that $H$ is a star Hopf algebra and that the action 
on $\CA$ is compatible with the star structures on both algebras:
\begin{equation}
\label{star-ac}
t \acts a^\ast = \left( (St)^\ast \acts a \right)^\ast.
\end{equation}

\begin{lemma}
The twisted algebra $H_\Psi$ is a star Hopf algebra 
provided that $\Psi^* = \Psi^{-1}$ (which for the twist
by the Cartan subalgebra translates to: $H_1^\ast \ts H_2^\ast 
= H_1 \ts H_2$).
\end{lemma}

This follows directly form (\ref{cotwi}).  

We shall prove that the action is compatible with the
star structure (\ref{star}), however to prove it  we shall 
rewrite (\ref{star}) in a more general way:
\begin{equation}
\ul{a}^\ast = \ul{ (U^{-1} \acts a)^\ast}, \label{star2}
\end{equation}
where $U$ is as defined in (\ref{utwi}). 

\begin{equation}
\begin{array}{l}
T \acts \ul{a}^\ast =  T \acts \ul{((U^{-1}  \acts  a)^\ast)} = \\
\phantom{xxx} =  \ul{ t \acts (U^{-1} \acts a)^\ast)} = \\
\phantom{xxx} =  \ul{ (St)^* U^{-1} \acts a)^\ast} 
\end{array}
\end{equation}
on the other hand:
\begin{equation}
\begin{array}{l}
T \acts \ul{a}^\ast =  ((ST)^\ast \acts \ul{a})^\ast = \\
\phantom{xxx} =  \left( (U (St) U^{-1})^\ast  \acts \ul{a} \right)^\ast = \\
\phantom{xxx} = \left( \ul{ U (St)^\ast U^{-1} \acts a} \right)^\ast =\\
\phantom{xxx} = \ul{ (U^{-1} U (st)^\ast U^{-1} \acts a)^\ast } = \\
\phantom{xxx} = \ul{ (St)^* U^{-1} \acts a)^\ast}.
\end{array}
\end{equation}
where we have used $U^{-1} = U^\ast$. This 
follows from $\Psi^* = \Psi^{-1}$ .

\subsection{The action of $U_\lambda({\mathfrak g})$}
For completeness we mention that one could always  
give the definition (\ref{action}) using the generators alone:
\begin{eqnarray}
&& {X_i^\pm} \acts \ul{T} = \ul{x_i^\pm \acts T}, \\
&& {H_i} \acts \ul{T} = \ul{h_i \acts T},
\end{eqnarray}
and then extend it to the whole of $U_\lambda({\mathfrak g})$
for every $T \in \CA_\lambda$.

To verify the compatibility of this action with the deformed
product in $\CA_\lambda$ clearly only action of $X_i^\pm$ 
must be checked, as the coproduct of $H_i$ remains 
not changed. Before we start  let us observe:
\begin{remark}
If $T \in \CA$ is homogeneous with degree $(n_1,n_2)$ 
then  $ x_i^\pm \acts T$ is also homogeneous with 
degree $(n_1 \pm a_{1i}, n_2 \pm a_{2i})$. This follows 
directly from (\ref{ug2}).
\end{remark}

\begin{proof}
Consider $X_i^\pm \acts (\ul{a} * \ul{b})$, where $a,b$ are 
homogeneous elements of $\CA$. On one hand:
\begin{equation}
X_i^\pm \acts (\ul{a} * \ul{b}) = 
\lambda^{n_1^a n_2^b} \left( \ul{ x_i^\pm \acts (ab) }\right),
\end{equation}
on the other hand, calculating it directly:
\begin{equation}
\begin{array}{rl}
X_i^\pm \acts (\ul{a} * \ul{b}) &= 
\left( X_i^\pm \acts \ul{a} \right) * \left( \lm{\mp \alpha_i H_2} \acts \ul{b} \right)
+ \left( \lm{\mp \beta_i H_1} \acts \ul{a} \right) * \left( X_i^\pm \acts \ul{b} \right) \\
&= \ul{(x_i^\pm \acts a)} * \ul{b} \lm{\mp \alpha_i n_2^b}  
+ \lm{\mp \beta_i n_1^a}   \ul{a} * \ul{(x_i^\pm \acts b)} \\
& = \lm{\mp \alpha_i n_2^b} \lm{ (n_1^a \pm \alpha_i) n_2^b }
\ul{(x_i^\pm \acts a) b} + \lm{\mp \beta_i n_1^a} \lm{ n_1^a (n_2^b \pm \beta_i)}
\ul{a(x_i^\pm \acts b)} \\
& = \lambda^{n_1^a n_2^b} \left( \ul{ x_i^\pm \acts (ab)} \right).
\end{array}
\end{equation}
\end{proof}
\section{Differential structures}

Differential structures on twisted Hopf algebras (and their
duals) has been extensively studied by Majid and Oeckl 
\cite{MaOe}, where the stability of bicovariant calculi under 
twisting has been shown. 
We briefly restate these results for the particular situation 
of the differential calculi over $\CA_\Psi$ invariant under
the action of $H_\Psi$.

Let $\Omega(\CA)$ be the differential algebra of forms over
the algebra $\CA$. We shall assume that the $\Omega(\CA)$ 
is invariant with respect to the action of $H$, i.e. the action 
of $H$ extends on $\Omega(\CA)$ and intertwines the 
exterior derivative:
\begin{equation}
H \ni t \acts (d \omega) = d( t \acts \omega),
\end{equation}

By using the same procedure as in the case of the deformation 
of $\CA$ we simply deform $\Omega(\CA)$ by modifying 
the product using (\ref{prod}). In particular, for two forms 
of homogeneous degree (defined similarly as for the elements 
of $\CA$) we have:
\begin{equation}
\ul{\omega} \wedge \ul{\rho} 
= \lm{n_\omega^1 n_\rho^2} \ul{(\omega \wedge \rho)},
\end{equation}
thus extending the quantization map to the differential algebra.

First we shall prove that on $\Omega(\CA_\Psi)$ there exists
an exterior derivative making it a differential algebra.

\begin{lemma}
The exterior derivative defined as:
\begin{equation} 
d \ul{\omega} = \ul{d \omega},  \label{diff-d}
\end{equation}
makes $\Omega(\CA_\Psi)$  a differential algebra.
\end{lemma}

We need to verify the Leibniz rule (it is obvious that $d^2 \equiv 0$). This follows, 
however, directly from the graded Leibniz rule of the undeformed differential
algebra and the relation (\ref{diff-d}).

\begin{lemma}
The deformed symmetry algebra acts on the quantized 
differential complex $\Omega(\CA_\Psi)$.
\end{lemma}

It remains only to verify that the action commutes 
with the external derivative.
$$
\begin{array}{l}
T \acts \ul{d \omega} = \ul{t \acts d \omega} = \\
\phantom{xxx} = \ul{d (t \acts \omega)}.
\end{array}
$$
and using (\ref{diff-d}) we see that the action 
of $T$ is well-defined.

Similarly one verifies that the star structure (if defined 
for the original differential complex) is preserved.
 
\section{Spectral triples}

Let us assume that there exists a spectral triple for the algebra
$\CA$, i.e. we have all the data : Hilbert space $\CH$ on which 
$\CA$ is represented as bounded operators, the Dirac operator $D$ (grading 
$\gamma$ in case of even dimensions and $J$ for real spectral triples) 
satisfying all the axioms as defined in \cite{Co2}. Additionally, we assume 
that the spectral triple is {\em symmetric} (as defined in \cite{PaSi}), so that 
there exists the action of an Hopf algebra $H$ on $\CA$, and that the crossproduct of $H$ and $\CA$ is represented on the Hilbert space. 
We call $H$ isometry if the representation of $H$ commutes with $D$.

Now we can state the main lemma:
\begin{lemma}
The deformation of the algebra $\CA$ by a twist $\Psi$ by a Cartan subalgebra 
allows for the representation of $\CA_\Psi$ and the crossproduct of $H_\Psi$ with 
$\CA_\Psi$ on the same Hilbert space. Moreover, with the Dirac operator $D$, as 
taken from the undeformed spectral triple we shall have a spectral triple for 
$\CA_\Psi$, which, moreover, will be invariant under the twisted Hopf algebra $H_\Psi$.
\end{lemma}

The sketch of the construction of the deformed spectral triple has been 
given in \cite{CoLa} (Theorem 6), here we only need to verify that $H_\Psi$
is the symmetry of algebra. 

We need to define the representation of $\CA_\Psi$ on the Hilbert space, which
gives rise to the representation of the crossproduct.

For $v \in \CH$, $\ul{a} \in \CA_\Psi$ we have:

\begin{equation}
\ul{a} v = \mu \left( \Psi^{-1} ( a \ts v) \right), \label{defrep}
\end{equation}

where $\mu$ denotes the representation of $\CA$, $\mu(a \ts v) = av$
Clearly, for $L \in H_\Psi$, $a \in \CA_\Psi$ and $v \in \CH$ we have:
\begin{equation}
L (\ul{a}v) = l \mu \left( \Psi^{-1}  \acts ( a \ts v) \right ) = \mu \left( (\cop l) \Psi^{-1}  \acts( a \ts v) \right ), 
\label{lhs1}
\end{equation}

where we have used the known form of the representation 
of the $H_\Psi$, which is the same as this of $H$   ($l$ 
is "undeformed" $L$). On the other hand:

$$ L (\ul{a}v) = \ul{\mu} \left( \Psi (\cop l) \Psi^{-1}  \acts (\ul{a} \ts v) \right),$$
where $\ul{\mu}$ is the representation of the deformed 
algebra $\mu(\ul{a} \ts v) =\ul{a}v$.
But using (\ref{defrep}) we obtain the same result
as in (\ref{lhs1}). \endproof

Since the representation of $H_\Psi$ is then the same as $\Psi$ 
and $D$ does not change as well, it follows at once that 
$H_\Psi$ will be the isometry of the deformed spectral triple.

\section{Conclusions}

As we have shown, the isospectral deformation has a Hopf 
algebra isometry, which is a twist of the classical Lie algebra
symmetry of the commutative space. Some physical models
based on twists have already been described in the literature, for 
instance, the twisted Lorenz group and Minkowski space, which 
has been discussed by Kulish and Mudrov in \cite{Kul}.

Since the algebraic form of the symmetries remains intact and 
only the coproduct changes one could expect modifications with 
respect to "undeformed physics" not on the one-particle level but 
only on the level of interactions. It remains to verified, especially 
for gauge theories, whether such deformations are physically 
plausible.

{\bf Acknowledgements:}
While finishing this paper we have learned of a similar work 
by Joseph Varilly \cite{Var} who developed the description
of the symmetries dual to the presented above. 

The author would also like to thank Michel Dubois-Violette, 
J.M.Gracia-Bond\'{\i}a, Gianni Landi,  John Madore, Mario Pachke 
and Harold Steinacker for helpful discussions on many topics.


\begin{thebibliography}{C}
\bibitem{CP} V.Chari, A.Pressley, {\em A Guide to Quantum Groups}, 
Cambridge University Press, 1994
\bibitem{Co} A.Connes, {\em Noncommutative geometry}, Academic Press 1994,
\bibitem{Co2} A.Connes, {\em Noncommutative geometry Year 2000}, arXiv:math.QA/0011193,
\bibitem{CoLa} A.Connes, G.Landi {\em Noncommutative manifolds, 
the instanton algebra and isospectral deformations}, arXiv:math.QA/0011194,
\bibitem{CMDV} A.Connes, M. Dubois-Violette, in preparation, 
\bibitem{GLS} A.Grossman, G.Loupias, E.M.Stein, {\em An algebra 
of pseudodifferential operators and quantum mechanics 
in phase space}, Ann.Inst.Fourier 18, 2, 343-368, (1969)
\bibitem{MDV} M. Dubois-Violette, {\em On the theory of quantum groups}, 
Lett.Math.Phys. 19, 121--126, (1990) 
\bibitem{Mor} E.Mourre, {\em Remarques sur le caractere algebraique du procede 
pseudo-differentiel et de certaines de ses extensions}, Ann.Inst.Henri 
Poincare, Phys.Theor. 53, 3, 259-182 (1990)
\bibitem{Kul} P.P.Kulish, A.I.Mudrov, {\em Twist-like geometries on a quantum 
Minkowski space.}, On the occasion of the 65-th birthday of Academician 
Lyudvig Dmitrievich Faddeev, Tr. Mat. Inst. Steklova 226, Mat. 
Fiz. Probl. Kvantovoi Teor. Polya, 97--111, (1999)
\bibitem{Maj} S.Majid, {\em Foundations of Quantum Group Theory}, 
Cambridge University Press, (1995)
\bibitem{MaOe} S.Majid, R.Oeckl, {\em Twisting of quantum differentials 
and the Planck scale Hopf algebra.}, Comm. Math. Phys. 205, 617--655,  (1999)
\bibitem{Pas} M.Paschke, {\em \"Uber nichtkommutative Geometrien, ihre Symmetrien 
und etwas Hochenergiephysik}, Ph.D. thesis, Mainz 2001, to appear, 
\bibitem{PaSi} M.Paschke, A.Sitarz, {\em The geometry of 
noncommutative symmetries},  Acta Phys.Pol. B31, No 11,  (2000)
\bibitem{Rie1} M.Rieffel, {\em Compact quantum groups associated with 
toral subgroups}, in: Representation Theory of group and Algebras, 
Contemp.Math. 145, AMS, Providence, R.I., (1993), 
\bibitem{Rie2} M.Rieffel, {\em Non-compact quantum groups associated with 
Abelian subgroups},  Comm.Math.Phys. 171, 181--201, (1995)
\bibitem{Rie3} M.Rieffel, {\em Deformation quantization for actions of 
$R^d$}, Memoirs AMS, vol 506, AMS, Providence, (1993)
\bibitem{Sit} A.Sitarz, {\em Rieffel's deformation quantization 
and isospectral deformations}, arXiv:math.QA/0102075, 
\bibitem{Var} J.C.V\'arilly, {\em Quantum symmetry groups of noncommutative 
spheres}, arXiv:math.QA/0102065, 
\end{thebibliography}
\end{document}